\documentclass{ita}
\usepackage{amsmath,amsfonts,amscd,amssymb,amsthm}

\def\Z{\mathbb Z}
\def\R{\mathbb R}
\def\Q{\mathbb Q}
\def\N{\mathbb N}
\def\A{\mathcal A}

\def\C{\mathcal C}
\def\L{\mathcal L}

\def\pf{\begin{proof}}
\def\pfk{\end{proof}}

\newtheorem{lem}{Lemma}[section]
\newtheorem{prop}[lem]{Proposition}
\newtheorem{coro}[lem]{Corollary}

\theoremstyle{definition}
\newtheorem{de}[lem]{Definition}
\newtheorem{pozn}[lem]{Remark}
\newtheorem{ex}[lem]{Example}

\begin{document}

\title{ On a class of infinite words with affine factor complexity}

\author{Julien Bernat}
\address{Institut de Math\'ematiques de Luminy  - CNRS UMR 6206\\ 163, avenue de Luminy,
case 907, 13288 Marseille cedex 09, France,
\email{bernat@iml.univ-mrs.fr}}
\author{Zuzana Mas\'akov\'a}
\address{Doppler Institute for Mathematical Physics and Applied
Mathematics \& Department of Mathematics, FNSPE, Czech Technical
University, Trojanova 13, 120 00 Praha 2, Czech Republic\\
\email{masakova@km1.fjfi.cvut.cz \&\ pelantova@km1.fjfi.cvut.cz}}
\author{Edita Pelantov\'a}
\sameaddress{2}



\begin{abstract}
In this article, we consider the factor complexity of a fixed
point of a primitive substitution canonically defined by a
$\beta$-numeration system. We provide a necessary and sufficient
condition on the R\'enyi expansion of 1 for having an affine
factor complexity map $\C(n)$, that is, such that $\C(n)=an+b$ for
any $n \in \N$.
\end{abstract}

\maketitle

\runningtitle{Affine complexity of beta-integers}

\runningauthors{Bernat, Mas\'akov\'a, and Pelantov\'a}


\section{Introduction}

Factor complexity is one of the basic properties which is studied
on infinite words $(u_n)_{n\in\N}$ over a finite alphabet
${\mathcal A}$. It is a function $\C:\N\to\N$, which counts the
number of factors of a given length which occur in an infinite
word. In other words, factor complexity expresses the measure of
irregularity in the word.

For eventually periodic words, the factor complexity is a function
bounded by a constant. As shown by Hedlund and Morse~\cite{morse},
an infinite word $(u_n)_{n\in\N}$ which is not eventually
periodic, i.e. is aperiodic, has factor complexity satisfying
$\C(n)\geq n+1$ for all $n\in\N$. Moreover, the language of the
factors of an infinite word is factorial, that is, one has
$\C(n+m) \leq \C(n) \C(m)$ for all $n,m \in \N$. It is therefore
obvious that not every function $\C$ can represent the factor
complexity of an infinite word. For an overview of necessary
conditions for a factor complexity function $\C$,
see~\cite{Ferenczi}.

Aperiodic words with minimal complexity $\C(n)=n+1$, for all
$n\in\N$, are called sturmian; their properties have been studied
by many authors, see~\cite{lothaire}. On the other hand, words
having maximal complexity satisfy $\C(n)=m^n$, where $m$ is the
cardinality of the alphabet. Under the term infinite words of low
factor complexity, one usually understands words for which $\C$ is
a sublinear function, i.e. there exist constants $a,b$ such that
$\C(n)\leq an+b$ for all $n\in\N$. A special subclass is formed by
infinite words with affine complexity, i.e.\ such that
$\C(n)=an+b$ for all $n\in\N$. Among the words with affine factor
complexity, one finds sturmian words, Arnoux-Rauzy words, words
coding generic interval exchange transformation, and others.

As shown by Queff\'elec~\cite{queffelec}, fixed points of a
primitive substitution have low factor complexity. Let us mention,
that relaxing the assumption of primitivity, the factor complexity
is bounded by a quadratic function, see~\cite{Pansiot}. The
determination of the factor complexity of a fixed point from the
prescription of the substitution is not a simple task.

\medskip

In this paper we consider canonical substitutions associated with
simple Parry numbers $\beta$. These are numbers whose R\'enyi
expansion of $1$ is finite, i.e.\ is of the form
$d_\beta(1)=t_1\cdots t_m$. The canonical substitution
corresponding to $\beta$ is a substitution over the alphabet
$\A=\{0,1,\dots,m-1\}$, given by
\begin{equation}\label{e:1}
\begin{array}{ccl}
\varphi(0)&=&0^{t_1}1\\
\varphi(1)&=&0^{t_2}2\\
&\vdots&\\
\varphi(m-2)&=&0^{t_{m-1}}(m-1)\\
\varphi(m-1)&=&0^{t_m}
\end{array}
\end{equation}

Since $t_1 \geq 1$ and $t_m \geq 1$, one easily checks that for
any letter $i$, $\varphi^{2m}(i)$ contains at least one occurence
of each letter, hence the substitution is primitive. Moreover, the
substitution $\varphi$ admits a unique fixed point, which is the
infinite word
$$
u_\beta:=\lim_{n\to\infty}\varphi^n(0)\,.
$$
In~\cite{FrMaPe}, the factor complexity of such fixed points is
determined for substitutions satisfying the condition 
$t_1>\max\{t_2,\dots,t_{m-1}\}$. In particular, one shows that
$$
(m-1)n+1\leq \C(n) \leq mn\,,\qquad \hbox{for all }\ n\geq 1\,.
$$
In the same paper it is shown that the word $u_\beta$ is
Arnoux-Rauzy, if and only if $t_m=1$ and $t_1=t_2=\cdots
=t_{m-1}$. In this case the factor complexity is obviously an
affine function.

\medskip

The aim of this article is the characterization of substitutions
of the form~\eqref{e:1}, for which the fixed point $u_\beta$ has
affine factor complexity. We will show

\begin{thrm}\label{t:1}
Let $\beta$ be a simple Parry number with the R\'enyi expansion of
unity $d_\beta(1)=t_1\cdots t_m$, and let $u_\beta$ be the fixed
point of the substitution~\eqref{e:1}. Then the factor complexity
of $u_\beta$ is an affine function if and only if the coefficients
$t_1,\dots,t_m$ satisfy
\begin{enumerate}
\item[1)] $t_m=1$
\item[2)] If there exists a non-empty word $w$ and $\alpha\in\Q$ such that
$t_1\cdots t_{m-1}=w^\alpha$, then $\alpha\in\N$.
\end{enumerate}
\end{thrm}

Let us mention that condition 2) of the above theorem means that
either $t_1\cdots t_{m-1}$ is equal to $w^k$ for $k\in\N$, $k\geq
2$, or no word can be both a proper prefix and a proper suffix of
$t_1\cdots t_{m-1}$. This formulation of condition 2) will be used
in the proof of the theorem.

Note that infinite words $u_\beta$ which are Arnoux-Rauzy, satisfy 
the condition $2)$ of the above theorem with $w = \lfloor \beta 
\rfloor$. Condition $2)$ is satisfied also by other words $u_\beta$, 
which are not Arnoux-Rauzy, but have the same complexity 
$\C(n)=(m-1)n+1$. These words illustrate the fact that Arnoux-Rauzy 
words of order $m\geq 3$ cannot be characterized by their complexity, 
as is the case for Arnoux-Rauzy words of order $m=2$, i.e.\ sturmian 
words.

In order to prove that conditions $1)$ and $2)$ of
Theorem~\ref{t:1} are sufficient for affine factor complexity, we
use purely the tools of combinatorics on words. For the opposite
implication, we use the geometric representation of the factors of
the word $u_\beta$ as coding of patterns occurring in the set of
$\beta$-integers, see section~\ref{sec:preli}.

\section{Preliminaries}\label{sec:preli}


\subsection{$\beta$-numeration}

In~\cite{renyi} the author introduces and studies the properties
of the positional number system with the base $\beta\in\R$,
$\beta>1$. For arbitrary real $x>0$, the $\beta$-expansion of $x$
can be found by the greedy algorithm, as follows. There exists a
unique $k\in\N$ such that $\beta^k\leq x<\beta^{k+1}$. Set
$x_k:=\lfloor x/\beta^k\rfloor$ and $r_k:=x-x_k\beta^k$. For each
$i<k$, set $x_i:=\lfloor\beta r_{i+1}\rfloor$ and $r_i:=\beta
r_{i+1}-x_i$. Obviously,
\begin{equation}\label{e:2}
x=x_k\beta^k + x_{k-1}\beta^{k-1} + x_{k-2}\beta^{k-2} + \cdots
\end{equation}
and $x_i\in\{0,1,\dots,\lceil \beta \rceil -1\}$. Note that the
elements of the sequence $(x_i)_{i \leq k}$ satisfy the relation
$x_i = \bigl\lfloor \beta T_{\beta}^{k-i}(x \beta^{-(k+1)})
\bigr\rfloor$, where the map $T_{\beta}$ is defined as:
\begin{equation}
T_{\beta}: [0,1] \rightarrow [0,1),\qquad T_\beta(x) = \beta x \!
\! \mod 1.
\end{equation}

For the expression of $x$ in the form of its
$\beta$-expansion~\eqref{e:2} we use the notation $x=x_k\cdots
x_0\bullet x_{-1}x_{-2}\cdots$, if $k\geq 0$, or $x=0\bullet
\underbrace{000\cdots 0}_{-k-1\ \hbox{\tiny times}}x_k
x_{k-1}\cdots$ if $k<0$. If the $\beta$-expansion ends in
infinitely many 0's, we omit them.

Numbers $x$ with vanishing $\beta$-fractional part, i.e.\ such
that $x_i=0$ for $i<0$ are called non-negative $\beta$-integers
and we denote them $x=x_k\cdots x_1x_0\bullet$. The set of
non-negative $\beta$-integers is denoted by $\Z_\beta^+$, and the
set of $\beta$-integers is defined as $\Z_{\beta} = \Z_{\beta}^+
\cup (-\Z_{\beta}^+) $.

Unlike the situation with integer base $\beta$, in case that
$\beta\notin\N$, there exist sequences $(x_i)_{i\leq k}$,
$x\in\{0,1,\cdots,\lceil\beta\rceil-1\}$ that are not the
$\beta$-expansion of some $x>0$. For the description of admissible
sequences of digits, one needs the so-called R\'enyi expansion of
1. For $\beta\in\R$, $\beta>1$, put $t_1:=\lfloor\beta\rfloor$ and
let $0\bullet t_2t_3t_4\cdots$ be the $\beta$-expansion of the
number $\beta-\lfloor\beta\rfloor$. Then the sequence
$d_\beta(1)=t_1t_2t_3\cdots$ is called the R\'enyi expansion of 1.
We have obviously,
$$
1=\sum_{i=1}^\infty \frac{t_i}{\beta^i} \quad\hbox{ and }\quad
t_i\in\{0,1,\dots,\lceil\beta\rceil-1\}.
$$
In order that a sequence $t_1t_2t_3\cdots$ of integers be the
R\'enyi expansion of $1$ for some base $\beta$, the so-called
Parry condition must be satisfied~\cite{parry},
\begin{equation}\label{e:3}
t_it_{i+1}t_{i+2}\cdots\prec t_1t_2t_3\cdots = d_\beta(1) 
\quad\hbox{for all } i\in\N,\ i\geq 2\,,
\end{equation}
where the symbol $\prec$ stands for `strictly lexicographically
smaller'. In the same paper~\cite{parry} it is shown that a finite
sequence of digits $x_kx_{k-1}\cdots x_1x_0$ over the alphabet
$\A=\{0,1,\dots,\lceil\beta\rceil-1\}$ is the $\beta$-expansion of
a $\beta$-integer if and only if
\begin{equation}\label{e:4}
x_ix_{i-1}\cdots x_0\prec d_\beta(1) \quad\hbox{ for all }
i\in\N,\ i\leq k\,.
\end{equation}

Using the R\'enyi expansion of $1$, one can even describe the
distances between consecutive $\beta$-integers on the real line.
If $\beta\in\N$, the $\beta$-integers are precisely the rational
integers, therefore the distance between consecutive
$\beta$-integers is always $1$. The situation is very different if
the base $\beta$ is not an integer. The distances between
consecutive $\beta$-integers are the elements of
$\bigl\{T_{\beta}^i(1) \mid {i \in \N} \bigr\}$,
see~\cite{thurston}. Note that, since $\Z_{\beta}$ is a discrete
set for any $\beta>1$, one may define the successor and
predecessor maps, respectively as
$$
{\rm pred}(x)=\max\{y \in \Z_{\beta}\mid y<x\}\quad\hbox{ and
}\quad {\rm succ}(x)=\min\{y \in \Z_{\beta}\mid y>x\}\,.
$$

\begin{ex}
Consider the base $\beta=\frac12(1+\sqrt5)$, i.e.\ $\beta$ is the
golden ratio. The number $\beta$ is a root of the equation
$x^2=x+1$, and its R\'enyi expansion is equal to
$d_{\beta}(1)=11$. The condition~\eqref{e:4} implies in this case
that $x_k\cdots x_0\in\{0,1\}^{k+1}$ is a $\beta$-expansion of a
$\beta$-integer, if and only if $x_ix_{i-1}\neq 11$ for all
$i=1,\dots,k$. The set of $\beta$-integers thus starts with the
numbers (written in their $\beta$-expansion)
$$
0\bullet\,,\quad 1\bullet\,,\quad 10\bullet\,,\quad
100\bullet\,,\quad 101\bullet\,,\quad 1000\bullet\,,\quad {\it
etc.}
$$
The distances between consecutive $\beta$-integers take only two
values, namely $T_{\beta}^0(1)=1$ and $T_{\beta}(1)=\beta-1$, see
Figure~\ref{f}.
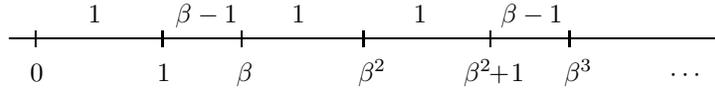
\begin{figure}
\begin{picture}(290,40)
 \put(10,20){\line(1,0){270}}
 \put(20,17){\line(0,1){6}}
 \put(18,4){0}
  \put(40,26){$1$}
 \put(68,17){\line(0,1){6}}
 \put(66,4){1}
  \put(73,26){$\beta-1$}
 \put(98,17){\line(0,1){6}}
 \put(96,4){$\beta$}
  \put(117,26){$1$}
 \put(144,17){\line(0,1){6}}
 \put(142,4){$\beta^2$}
  \put(163,26){$1$}
  \put(192,17){\line(0,1){6}}
 \put(182,4){$\beta^2\!\!+\!1$}
  \put(196,26){$\beta-1$}
  \put(222,17){\line(0,1){6}}
 \put(220,4){$\beta^3$}
  \put(260,4){$\cdots$}
\end{picture}
\caption{The set of $\beta$-integers for $\beta=\frac12(1+\sqrt5)$
drawn on the real line.} \label{f}
\end{figure}
\end{ex}

Numbers $\beta$, for which the R\'enyi expansion $d_\beta(1)$ is
eventually periodic, are called Parry numbers. In this case the
number of values for distances between consecutive
$\beta$-integers is finite. In other words, Parry numbers are
numbers for which the elements of the sequence
$(T_{\beta}^i(1))_{i \in \N}$ take finitely many distinct values.
If moreover $d_\beta(1) = t_1 t_2 \cdots t_m $, $t_m\neq 0$, then
$\beta$ is called a simple Parry number; one has $T_{\beta}^i(1) =
0 \bullet t_{i+1} \cdots t_m $ for any $i \in \{1,\ldots,m-1\}$.
Associating to these distances letters in the alphabet
$\A=\{0,1,\dots,m-1\}$ in a natural way, $T_{\beta}^i(1)\mapsto
i$, the distances between consecutive $\beta$-integers can be
coded by an infinite word over $\A$. This infinite word is the
fixed point of the substitution $\varphi$ defined by~\eqref{e:1}.

\begin{ex}
The infinite word $u_\beta$ for $\beta=\frac12(1+\sqrt5)$ starts
by
$$
u_\beta=01001\cdots
$$
\end{ex}

In~\cite{fabre} it is shown that the infinite word $u_\beta$ is
a fixed point of a canonical substitution~\eqref{e:1} associated
to $\beta$. Note that a canonical substitution can be associated
also to a non-simple Parry number $\beta$, see~\cite{fabre}.
For more details about the properties of $\beta$-numeration we
refer to~\cite{lothaire}.

\subsection{Combinatorics on words}

Let $\A=\{0,1,\dots,m-1\}$ be a finite alphabet. A finite
concatenation $w=w_0w_1\cdots w_{n-1}$ of the letters is called a
word, its length $n$ is denoted by $|w|$. The set of finite words
over an alphabet $\A$ together with the empty word $\varepsilon$
and the concatenation operation forms a free monoid, denoted by
$\A^*$.

The sequence $u=(u_i)_{i\in\N}$ of the letters in the alphabet
$\A$ is called an infinite word. A word $w$ is a factor of a word
$u$ (finite or infinite), if there exist words $w^{(1)}$ and
$w^{(2)}$ such that $u=w^{(1)}ww^{(2)}$. If $w^{(1)}$ is an empty
word, then $w$ is a prefix of $u$. If $w^{(2)}=\varepsilon$, then
$w$ is a suffix of $u$. The set of all factors of an infinite word
$u$ is called the language of $u$ and denoted by $\L(u)$. The set
of all factors of $u$ of length $n$ is denoted by $\L_n(u)$.
Obviously $\L(u)=\bigcup_{n\in\N}\L_n(u)$. The cardinality of the
set $\L_n(u)$ is the factor complexity. Formally, we have the
function $\C:\N\to\N$, given by
$$
\C(n):=\#\L_n(u)\,.
$$
Note that any language which consists of the set of factors of an
infinite words is extendable, that is, every factor $w_0\cdots
w_{n-1}$ of length $n$ can be extended in at least one way to a
factor $w_0\cdots w_{n-1}w_n$ of length $n+1$. Hence the factor
complexity is a non-decreasing function. The set of letters, by
which it is possible to extend a factor $w$ to the right is called
the right extension of $w$,
$$
{\rm Rext}(w) = \{a\in\A \mid wa\in\L(u)\}\,.
$$
The increment of complexity can be calculated using the number of
right extensions of all factors of length $n$,
$$
\Delta\C(n) = \C(n+1)-\C(n) = \sum_{w\in\L_n(u)} \bigl(\#{\rm
Rext}(w) -1\bigr)\,.
$$
A factor $w$, for which $\#{\rm Rext}(w)\geq 2$ is called a right
special factor. Only such factors are important for the
determination of the first difference of factor complexity.

In this paper we study recurrent words. These are infinite words,
in which every factor appears at least twice. Factors of a
recurrent word can be extended in at least one way to the left,
and so  all the above considerations can be analogically stated.
In particular, we have
\begin{equation}\label{e:5}
\Delta\C(n) =  \sum_{w\in\L_n(u)} \bigl(\#{\rm Lext}(w)
-1\bigr)\,,
\end{equation}
where ${\rm Lext}(w) = \{a\in\A \mid aw\in\L(u)\}$. Factors with
$\#{\rm Lext}(w)\geq 2$ are called left special. Factors which are
both right special and left special are called bispecial.

A morphism on the free monoid $\A^*$ is a mapping
$\varphi:\A^*\to\A^*$ satisfying
$\varphi(wv)=\varphi(w)\varphi(v)$ for every pair $w,v\in\A^*$. It
is obvious that a morphism is uniquely determined by images of all
letters $a\in\A$. The action of a morphism can be naturally
extended to infinite words $(u_n)_{n\in\N}$ as
$$
\varphi(u)=\varphi(u_0u_1u_2\cdots)
:=\varphi(u_0)\varphi(u_1)\varphi(u_2)\cdots
$$
If moreover $\varphi(a) \neq \varepsilon$ for all $a\in\A$, and
there exist $a_0$ which is a proper prefix of $\varphi(a_0)$, then
the morphism $\varphi$ is called a substitution. An infinite word
$u$ satisfying $u = \varphi(u)$ is a fixed point of the
substitution $\varphi$. Obviously, a substitution has at least one
fixed point, namely $\lim_{n\to\infty}\varphi^n(a_0)$. A
substitution is called primitive, if there exists $k \in \N$ such
that, for every pair of letters $a,b \in \A$, the letter $a$
appears in the word $\varphi^k(b)$. It is known~\cite{Durand} that
a fixed point of a primitive substitution is a linearly recurrent
word, which implies that the distances between consecutive
occurrences of a given factor are bounded.

\section{Affine factor complexity of infinite words $u_\beta$}

Our aim is to describe the substitutions of the form~\eqref{e:1}
whose fixed points
$$
u_\beta = \underbrace{0^{t_1}1\ 0^{t_1}1\ \cdots \ 0^{t_1}1}_{t_1
\hbox{ \small times}}\ 0^{t_2}2 \cdots
$$
have affine factor complexity, i.e.\ the first difference
$\Delta\C(n)$ is constant. For the determination of $\Delta\C(n)$
we use the left special factors of $u_\beta$. In~\cite{FrMaPe} it
is shown that every prefix $w$ of the infinite word $u_\beta$ is a
left special factor and its left extension is ${\rm
Lext}(w)=\A=\{0,1,\dots,m-1\}$. Therefore using~\eqref{e:5} we
have $\Delta\C(n)\geq m-1$ for every $n\in\N$.

Infinite words whose every prefix is a left special factor are
called left special branch~\cite{cassaigne}. As we have mentioned,
$u_\beta$ is a left special branch of itself. In~\cite{FrMaPe} it
is moreover shown that $u_\beta$ has no other left special branch.

For the description of left special factors of another type (i.e.\
which are not prefixes of a left special branch) we use a lemma
from~\cite{FrMaPe}.

\begin{de}
Let $d_\beta(1)=t_1t_2\cdots t_m$. For $2\leq k\leq m$ we denote
$$
j_k:=\min\{i\in\N\mid 1\leq i\leq k-1,\ t_{k-i}\neq 0\}.
$$
\end{de}

Note that an index $j_k$ always exists, because $t_1>0$.

\begin{lem}\label{l:1}
All factors of $u_\beta$ of the form $X0^rY$, where $X,Y$ are
non-zero letters and $r\in\N$, are the following,
$$
\begin{array}{cl}
j_k0^{t_k}k\,,&\hbox{ for } k=2,3,\dots,m-1\,,\\[2mm]
k\,0^{t_1}1\,,&\hbox{ for } k=1,2,\dots,m-1\,,\\[2mm]
j_m0^{t_1+1}1\,.
\end{array}
$$
\end{lem}

This lemma is exactly Lemma 4.5 in ~\cite{FrMaPe}.

\begin{pozn}\label{R}\
\begin{itemize}
\item[(i)] Since $j_k\leq k-1$, the only factors of the form $(m-1)0^rY$
are, according to Lemma~\ref{l:1}, the factors $(m-1)0^{t_1}1$ and
possibly $j_m0^{t_1+t_m}1$. In any case, the letter $(m-1)$ is
always succeeded by the letter $0$.

\item[(ii)] Recall that for parameters $t_1,\dots,t_m$ of the
substitution it holds that $t_m\geq 1$, and from the Parry
condition $t_1\geq t_i$ for all $i=2,\dots,m$.
\end{itemize}
\end{pozn}

\begin{coro}\label{c:2}
Every left special factor $w$ with $|w|\leq t_1$ is a prefix of
$u_\beta$.
\end{coro}

\pf We prove the statement by contradiction. Let $w$ be a left
special factor satisfying $|w|\leq t_1$, and suppose that $w$ is
not a prefix of $u_\beta$. Since $u_\beta$ has a prefix $0^{t_1}$,
necessarily $w$ is of the form $0^rY$ for some $Y\neq0$, $r\in\N$,
$r<t_1$. If $Y=1$, then from Lemma~\ref{l:1} we know that $w$ has
a unique left extension, namely 0, and thus cannot be a left
special factor. If $Y>1$, then again, $w$ has a unique left
extension, namely 0, if $r<t_Y$, or $j_k$, if $r=t_Y$.
\pfk

\begin{prop}\label{p:3}
Let $\beta$ be a simple Parry number. The infinite word $u_\beta$
has affine factor complexity if and only if every left special
factor is a prefix of $u_\beta$.
\end{prop}

\pf Since every prefix $w$ of $u_\beta$ satisfies $\#{\rm
Lext}(w)=m$, Corollary~\ref{c:2} implies that $\Delta\C(n)=m-1$
for all $n\leq t_1$. If $u_\beta$ has affine factor complexity,
then $\Delta\C(n)=m-1$ for all $n\in\N$, and so no left special
factors other than prefixes of $u_\beta$ can exist. The opposite
implication is obvious. \pfk

\begin{coro}\label{c:4}
If $u_\beta$ has affine factor complexity, then $t_m=1$.
\end{coro}

\pf Suppose that $t_m\geq 2$. Then according to Lemma~\ref{l:1},
the word $0^{t_1+t_m-1}$ is a left special factor, because it has
two distinct left extensions, namely 0 and $j_m$. In the same
time, $0^{t_1+t_m-1}$ is not a prefix of $u_\beta$.
Proposition~\ref{p:3} implies that the factor complexity of
$u_\beta$ is not an affine function. \pfk

In~\cite{FrMaPe} it is shown that under the conditions
$$
(a)\quad t_m=1 \qquad\qquad (b)\quad t_1=t_2=\cdots = t_{m-1}
\quad \hbox{ or }\quad t_1>\max\{t_2,\dots,t_{m-1}\}\,,
$$
%
%
%
the factor complexity of $u_\beta$ is affine. Note that
the condition $(b)$ is a very special case of condition 2) of
Theorem~\ref{t:1}, whose proof is the aim of this paper.

\begin{de}  A left special  factor $w$ of an infinite word $u$ is
called maximal if for any letter   \ $a\in\A$  \  the word $wa$ is
not a left special factor of $u$.
\end{de}

If $t_m\geq 2$, then $0^{t_1+t_m-1}$ is maximal, since extending
it to the right using Lemma~\ref{l:1}, we do not obtain a left
special factor. Let us mention that if $w$ is a maximal left
special factor, then it is a bispecial factor: Since $w$ is left
special, there exist $X_1,X_2\in\A$ such that $X_1w,
X_2w\in\L(u_\beta)$. Every factor of $u_\beta$ can be extended in
at least one way to the right, and thus we can find $Y_1,Y_2\in\A$
so that $X_1wY_1$ and $X_2wY_2$ belong to $\L(u_\beta)$. Since $w$
is a maximal left special factor, we have $Y_1\neq Y_2$. This
however means that $w$ is a right special factor.

Every left special factor $w$ is either maximal or it can be
extended by a letter $a\in\A$ such that $wa$ is again a left
special factor. Since the only infinite left special branch of
$u_\beta$ is $u_\beta$ itself, every left special factor which is
not prefix of $u_\beta$ is a prefix of a maximal left special
factor. Proposition~\ref{p:3} therefore implies the following
Corollary.

\begin{coro}\label{c:5}
The infinite word $u_\beta$ has affine factor complexity if and
only if $u_\beta$ has no maximal left special factor.
\end{coro}

\subsection{Sufficient condition for affine factor complexity of $u_\beta$}

In the previous part we have derived that $u_\beta$ can have
affine factor complexity only if $t_m=1$. Therefore we shall
consider only simple Parry numbers with the R\'enyi expansion
$$
d_\beta(1)=t_1t_2\cdots t_{m-1}1
$$
and study the substitution
\begin{equation}\label{e:6}
\begin{array}{ccl}
\varphi(0)&=&0^{t_1}1\\
\varphi(1)&=&0^{t_2}2\\
&\vdots&\\
\varphi(m-2)&=&0^{t_{m-1}}(m-1)\\
\varphi(m-1)&=&0
\end{array}
\end{equation}

In agreement with Corollary~\ref{c:5}, the study of conditions
under which the factor complexity is an affine function, resumes
into the study of existence of maximal left special factors in the
language of $u_\beta$. Lemma~\ref{l:1} under the condition $t_m=1$
states that the longest factor containing only zero letters is
$0^{t_1+1}$, and this factor has a unique extension to the left
and to the right. Therefore a left special factor of the form
$0^r$ satisfies $r\leq t_1$, and hence it is a prefix of the infinite
left special branch $u_\beta$.

We have thus shown the following simple observations.

\begin{lem}\label{l:6}
Any maximal left special factor contains at least one non-zero
letter.
\end{lem}

From the form of the substitution~\eqref{e:6} one can deduce the
structure of left special factors.

\begin{lem}\label{l:7}
If $w\in  \L(u_\beta)$ is a left special factor (not necessary maximal) then
$$
w=
\left\{%
\begin{array}{ll}
    0^r, & \hbox{for some}\ \ r\in \N ,\ r\leq t_1 ; \\
    \varphi(v)0^s, & \hbox{for some left special factor} \ v \  \hbox{and }  \ s\in \N .\\
\end{array}%
\right.
$$
\end{lem}

\begin{lem}\label{l:8} Let $w\in  \L(u_\beta)$.
\begin{enumerate}
\item If $w$ is a left special factor then $\varphi(w)$ is a left special factor
with the same number of left extensions;

\item  If $w$ is a maximal left special factor then there exists $q\in \N, q\leq
t_1$ such that $\varphi(w)0^q$ is a maximal left special factor.
\end{enumerate}
\end{lem}

The statement (2) of Lemma~\ref{l:8} says that if there exists one
maximal left special factor, then there exists an entire sequence
of them.

\begin{de}\label{initial}
A maximal left special factor $w$ is called non-initial if there
exists a maximal left special factor $v$ and an integer $q\in \N$
such that $w =\varphi(v)0^q$. A maximal left special factor which
is not non-initial is called initial maximal left special factor.
\end{de}

If $\L(u_\beta)$ contains a maximal left special factor, then it
contains an initial maximal left special factor as well. In order
to describe initial maximal left special factors, we introduce the
notion of trident.

\begin{de}\label{trident}
A factor $w \in \L(u_\beta)$ is called a trident if there exists
letters $X,Y,Z \in \A$ such that
\begin{enumerate}
\item $wX$ is a left special factor;
\item $wY$ and $wZ$ are not left special factors;
\item the unique left
 extensions of $wY$ and $wZ$ are distinct.
\end{enumerate}
The letter $X$ is called the rooted tooth, the letters $Y$ and $Z$
are called non-rooted teeth of the trident $w$.
\end{de}

Clearly, the teeth $X,Y,Z$ are different.

\begin{pozn}\label{R2}
If $0^r$ is a trident, then the rooted tooth $X=0$ or $1$. This
fact follows from Lemma~\ref{l:1}, since $0^rX$ is a left special
factor only if $X\leq 1$.
\end{pozn}

\begin{lem}\label{l:9}
Let $w$ be a trident containing a non-zero letter with rooted
tooth $X$ and non-rooted teeth $Y$, $Z$.

\begin{itemize}
\item[(i)] If $X=0$ then $t_Y=t_Z$.

\item[(ii)] If $X\neq0$ then there exists an integer $s\in\N$ and
a trident $\hat{w}$
with rooted tooth $\hat{X}\neq m-1$ and non-rooted teeth
$\hat{Y}$, $\hat{Z}$, such that

{\rm \ (a)} $w=\varphi(\hat{w})0^{s}$,

{\rm \ (b)} $A=\hat{A}+1$ for every non-zero tooth $A$ of the trident $w$,

{\rm \ (c)} $s=t_A$ for every non-zero tooth $A$ of the trident $w$,

{\rm \ (d)} if $A= 0$ is a non-rooted tooth of $w$, then $\hat{A}=m-1$ or $t_{\hat{A}+1}>t_X=t_{\hat{X}+1}$.
\end{itemize}
\end{lem}

\pf From the definition of a trident, it follows that $w$ is a
left special factor. According to Lemma~\ref{l:7}, there exist a
left special factor $\hat{w}$ and $s\in\N$ such that
$w=\varphi(\hat{w})0^s$.

(i) Let $X=0$. Since $wY=\varphi(\hat{w})0^sY$ and
$wZ=\varphi(\hat{w})0^sZ$ are factors of $u_\beta$, and $Y,Z\neq
X=0$, it follows that $s=t_Y=t_Z$.

(ii) Let $X\neq 0$. Since
$wX=\varphi(\hat{w})0^sX=\varphi\bigl(\hat{w}(X-1)\bigr)$ is a
left special factor, also $\hat{w}(X-1)$ is a left special factor
and $s=t_X$. As teeth $Y, Z$ are distinct, at least one of them is
non-zero, say $Y\neq 0$.  Since
$wY=\varphi(\hat{w})0^sY=\varphi\bigl(\hat{w}(Y-1)\bigr)$ is not a
left special factor, due to Lemma~\ref{l:8}, $\hat{w}(Y-1)$ is
also not a left special factor and $t_X=t_Y$.

If moreover $Z\neq0$, we have analogically $t_X=t_Z$ and
$\hat{w}(Z-1)$ is not a left special factor. As factors
$\varphi\bigl(\hat{w}(Y-1)\bigr)$ and
$\varphi\bigl(\hat{w}(Z-1)\bigr)$ have different left extensions,
also factors $\hat{w}(Y-1)$ and $\hat{w}(Z-1)$ have different left
extensions, and therefore $\hat{w}$ is a trident with teeth $X-1$,
$Y-1$, $Z-1$.

Suppose now that $Z=0$. Since $\varphi\bigl(\hat{w}(X-1)\bigr)$,
$\varphi\bigl(\hat{w}(Y-1)\bigr)$ and $\varphi(\hat{w})0^{t_X}0$
are factors of equal length, there must exist a letter
$\hat{Z}\neq X-1,Y-1$ such that $\hat{w}\hat{Z}\in\L(u_{\beta})$
and $\hat{w}\hat{Z}$ has a unique left extension. If $\hat{Z}\neq
m-1$, then $w0=\varphi(\hat{w})0^{t_X}0$ is a proper prefix of
$\varphi(\hat{w}\hat{Z})=\varphi(\hat{w})(\hat{Z}+1)$, and hence
$t_{\hat{Z}+1}\geq t_X+1$. \pfk

\begin{coro}\label{R3}
If $w$ is a trident with rooted tooth $X=1$ and $Y\neq 0$ is
a non-rooted tooth, then $t_Y=t_1$.
\end{coro}

Tridents play important role for existence of maximal left special
factors.

\begin{prop}\label{p:10}
Let $v$ be an initial maximal left special factor. Then there
exists a trident $w$ with rooted tooth $X$ and an integer $s \in \N$ such that
\begin{equation}\label{e:7}
\begin{aligned}
v=\varphi(w)0^s,&\quad X\neq 0,m-1, \quad  \\
\hbox{and}\ & \quad t_{X+1} < s=\min\{ t_{A+1} \mid A \hbox{ is a non-rooted tooth of } w,\ A\neq m-1\}\,.
\end{aligned}
\end{equation}
\end{prop}

\pf Let $v$ be an initial maximal left special factor, and let
$y',z'\in\A$ be its distinct left extensions. Denote by $Y'$ the
right extension of $y'v$ and by $Z'$ the right extension of $z'v$.
Since $v$ is a maximal left special factor, necessarily $Y'\neq
Z'$. According to Lemmas~\ref{l:6} and~\ref{l:8}, we have
$v=\varphi(w)0^s$ for some left special factor $w$ and some
$s\in\N$, $s\leq t_1$. Since $\varphi(w)0^sY'$ and
$\varphi(w)0^sZ'$ belong to the language $\L(u_\beta)$, there
exist distinct letters $Y,Z$ such that $wY$, $wZ\in\L(u_\beta)$,
and $wY$, $wZ$ have unique left extensions.

Since $v=\varphi(w)0^s$ is an initial maximal left special factor,
the left special factor $w$ is not maximal, and thus there exists
a letter $X$ such that $wX$ is a left special factor. This shows
that the factor $w$ is a trident with rooted tooth $X$ and
non-rooted teeth $Y,Z$.

Let us now show that $X\neq0,m-1$. Suppose that $X=0$. Then using
Lemma~\ref{l:8}, the factor $\varphi(wX)=\varphi(w)0^{t_1}1$ is
left special. Since $v=\varphi(w)0^s$ and $s\leq t_1$, it
implies that $v$ is a prefix of a left special factor
$\varphi(w)0^{t_1}1$, which is a contradiction with maximality of
$v$.

Suppose now that $X=m-1$. Then using (i) of Remark~\ref{R}, the factor $w(m-1)0$
is left special, and thus $\varphi(w)0^{t_1+1}$ is also a left
special factor. Again, we obtain a contradiction with the
maximality of $v$, since $v$ is then a proper prefix of another
left special factor.

The same reason leads us to the fact that $s>t_{X+1}$, because
otherwise $v=\varphi(w)0^s$ is a proper prefix of the left special
factor $\varphi(w)\varphi(X)=\varphi(X)0^{t_{X+1}}(X+1)$, where we
use that $X\neq m-1$.

It remains to determine the value of $s$. Since at least one of
the letters $Y',Z'$ is non-zero, say $Y'\neq0$, we have
$vY'=\varphi(w)0^sY'=\varphi(wY)$, and thus $Y'=Y+1$, $s=t_{Y+1}\leq t_1$
and $Y\neq m-1$. If moreover $Z'\neq0$, we have by the same
arguments that $s=t_{Z+1}=t_{Y+1}$, and $Z\neq m-1$. If $Z'=0$,
then either $Z=m-1$ or $Z\neq m-1$ and
$t_{Z+1}>s=t_{Y+1}$.
 \pfk

We are now in position to prove that condition 2) of
Theorem~\ref{t:1} is sufficient for $u_\beta$ having affine factor
complexity.

\begin{prop}\label{p:11}
Let $u_\beta$ be the infinite word associated to the Parry number
$\beta$ with $d_\beta(1)=t_1\cdots t_{m-1}1$. If $u_\beta$ does
not have affine factor complexity, then

\begin{enumerate}
\item[1)] there exists a non-empty word which is both a proper prefix
and a proper suffix of the word $t_1\cdots t_{m-1}$;

\item[2)] for every $k\in\N$, $k\geq 2$, and every word $w$ it holds
that $w^k\neq t_1\cdots t_{m-1}$.
\end{enumerate}
\end{prop}

\pf If the factor complexity of $u_\beta$ is not an affine
function, then there exists an initial maximal left special factor
$v$. According to Proposition~\ref{p:10}, there exists an integer
$s$ and a trident $w$ with rooted tooth $X$ and non-rooted teeth
$Y,Z$ satisfying conditions~\eqref{e:7}. Denote $l=X$.
Relations~\eqref{e:7} imply that $1\leq l< m-1$. We want to
construct $l$ tridents $w^{(1)}$, $w^{(2)}$, \dots, $w^{(l)}$ with
triples of teeth $(1,Y_1,Z_1)$, $(2,Y_2,Z_2)$, \dots,
$(l,Y_l,Z_l)$, and integers $s_1,s_2,\dots,s_l$ such that
$w^{(i)}$, $w^{(i+1)}$ and $s_{i+1}$ have properties of tridents
$\hat{w}$, $w$ and the integer $s$ from Lemma~\ref{l:9} for all
$i=1,\dots,l-1$, and $Y_l=Y$ and $Z_l=Z$. If $l=1$, this role is
played obviously by the trident $w$, its triple of teeth $(1,Y,Z)$
and the integer $s$. If $l\geq 2$, then according to
Remark~\ref{R2}, the trident $w$ contains a non-zero letter and
satisfies assumptions of Lemma~\ref{l:9}, which implies the
existence of the sequence of tridents $w^{(1)}$, $w^{(2)}$, \dots,
$w^{(l)}$ with triples of teeth $(1,Y_1,Z_1)$, $(2,Y_2,Z_2)$,
\dots, $(l,Y_l,Z_l)$, and integers $s_1,s_2,\dots,s_l$ with
required properties.

According to Corollary~\ref{R3}, we have
\begin{equation}\label{e2}
s_1=t_1.
\end{equation}
Since the rooted teeth $X_1=1$, $X_2=2$, \dots, $X_l=l$ are non-zero,
(c) of Lemma~\ref{l:9} implies
\begin{equation}\label{e1}
s_2=t_2,\quad s_3=t_3,\quad \dots,\quad s_l=t_l\,.
\end{equation}
Using Proposition~\ref{p:10} we obtain
\begin{equation}\label{e3}
t_{l+1} < s:=\min\{ t_{A+1} \mid A \hbox{ is a non-rooted tooth of } w^{(l)},\ A\neq m-1\}\,.
\end{equation}

Lemma~\ref{l:9} implies that the sequence $Y_1$, $Y_2$, \dots, $Y_l$ is formed by consecutive
integers separated by blocks of 0's. More precisely, for any $i=1,\dots,l-1$, we have
\begin{equation}\label{e4}
Y_{i+1}=\left\{\begin{array}{clcl}
Y_{i}+1 & \hbox{ if }\ \ Y_i<m-1 &\hbox{and} &t_{Y_i+1}=s_{i+1}\,,\\[2mm]
0       & \hbox{ if }\ \ Y_i=m-1 &\hbox{or}  &t_{Y_i+1}>s_{i+1}\,.
\end{array}\right.
\end{equation}
The same rule is valid for the sequence $Z_1$, \dots, $Z_l$.

Since non-rooted teeth $Y_1,Z_1$ are distinct, we can without loss
of generality assume that $Y_1\geq 2$.

In order to show the statement 1) of the proposition, denote by $k\leq l$ the maximal
index such that the sequence $Y_1$, \dots, $Y_k$ is formed by consecutive non-zero integers, i.e.\
$$
Y_1\ Y_2\ \cdots \ Y_k \ = \ j \ (j+1) \ \cdots \ (j+k-1)\qquad \hbox{ for some }\ j\in\N,\ j\geq 2\,.
$$
This however means, using~\eqref{e4},~\eqref{e1} and Corollary~\ref{R3} that
\begin{equation}\label{e5}
t_j=t_1,\quad t_{j+1}=t_2,\quad \dots,\quad t_{j+k-1}=t_k\,.
\end{equation}

We now show that the non-rooted tooth $Y_k=(j+k-1)$ is equal to $(m-1)$, which together with~\eqref{e5}
results in the statement (1) of the proposition.
For the contradiction, assume that $Y_k=(j+k-1)<m-1$. Let us distinguish two cases according to whether $k<l$
or $k=l$.
If $k<l$, then from the definition of $k$ it follows that $Y_{k+1}=0$, which, due to~\eqref{e4},
can happen only if
\begin{equation}\label{e6}
t_{Y_k+1}=t_{j+k} \ > \ s_{k+1}=t_{k+1}\,.
\end{equation}
If $k=l$, then~\eqref{e3} implies
\begin{equation}\label{e7}
t_{k+1}< s \leq t_{Y_{k}+1}=t_{j+k}\,.
\end{equation}
In any case,~\eqref{e5} together with~\eqref{e6}, or~\eqref{e7} gives
$$
t_{j}t_{j+1}\cdots t_{j+k} \succ t_1t_2\cdots t_{k+1}\,,
$$
which contradicts the Parry condition~\eqref{e:3}.

Besides the validity of the statement (1) of the proposition, we have thus
proved that the sequence $Y_1,\dots, Y_l$ contains at least one letter $m-1$.

\medskip
In order to show the statement 2) of proposition, denote by $p$
the shortest non-empty word which is both a proper prefix and a
proper suffix of the word $t_1\cdots t_{m-1}$. It is obvious that
$p$ is not a power of a shorter word.

We show the statement (2) by contradiction. Assume that there exists a word
$w$ such that $w^k=t_1\cdots t_{m-1}$ for some $k\geq 2$, $k\in\N$. First
we claim that such an assumption implies that $t_1\cdots t_{m-1}=p^n$ for
some $n\in\N$, $n\geq 2$. Since $w$ is a prefix and a suffix of $t_1\cdots t_{m-1}$,
we must have $|w|\geq p$. If $|w|=|p|$, the claim is valid. If $|w|>|p|$,
then $p$ is a proper prefix and a proper suffix of $w$.
Moreover, the prefix $p$ and the suffix $p$ do not overlap in the
word $w$, since otherwise the overlap would be a proper prefix and
a proper suffix of $t_1\cdots t_{m-1}$ shorter than $p$, which
contradicts the minimality of $p$. The condition $|w|>|p|$ thus implies
that $w=pw'p$ for some (possibly empty) word $w'$. If $w'=\varepsilon$,
the claim is valid.
In the opposite case, the word $t_1t_2\cdots t_{m-1}$ has the
prefix $pw'ppw'p$. The Parry condition for $d_\beta(1)$ implies
that $w'p\preceq pw'$, and $ppw'\preceq pw'p$ which then implies
$pw'\preceq w'p$, and therefore $pw'=w'p$. It is known that if two
words commute, then they are powers of the same word. Since $p$
itself is not a power, we must have $w'=p^j$ for some $j\in\N$, as
we wanted to show.

Let now $t_1t_2\cdots t_{m-1}=p^n$ for some $n\in\N$, $n\geq 2$.
Denote $s=|p|$. Obviously $m-1=ns$. If $s=1$, then
$d_\beta(1)=t_1t_1\cdots t_11$, and in that case $u_\beta$ is an
Arnoux-Rauzy word, for which it is known that the factor
complexity is an affine function. Thus $s\geq 2$.

Let us come back to the sequence of tridents and the triples of
their teeth, $(1,Y_1,Z_1)$, $(2,Y_2,Z_2)$, \dots, $(l,Y_l,Z_l)$.
We already know that one of the letters $Y_1$, \dots, $Y_l$ is
equal to $m-1$. Denote by $q$ the maximal index, such that $Y_q$
or $Z_q$ is equal to $m-1=ns$. Since the role of $Y_q$ and $Z_q$
is symmetric, without loss of generality we can assume that the
last $m-1$ occurred was $Y_q=m-1$. We will show that both the
corresponding rooted tooth $q$ and the other non-rooted tooth
$Z_q$ are multiples of $s$.

For a contradiction, suppose that $q=as+b$, where $1\leq b<s$.
According to Lemma~\ref{l:9}, we have
$$
t_q=t_{Y_q}=t_{m-1}, \quad t_{q-1}=t_{m-2},\quad \cdots,\quad
t_{q-b+1}=t_{m-1-b+1}\,.
$$
Since the word $p$ of the length $s$ is the period of $t_1\cdots t_{m-1}$, we have
$$
t_q=t_{as+b}=t_b=t_{m-1},\quad t_{b-1}=t_{m-2},\quad \cdots,\quad
t_{q-b+1}=t_1=t_{m-b}\,,
$$
and therefore $t_1\cdots t_b$ is both a prefix and a suffix of
$t_1\cdots t_{m-1}$, shorter than $p$, which contradicts the
choice of $p$. In the same way, one can show that the non-rooted
tooth $Z_q$ is a multiple of $s$, say $Z_q=cs$ for some $c\in\N$.

Since for the sequence of letters $Z_1$, \dots, $Z_l$ one can derive a rule
analogous to~\eqref{e4}, we obtain from the periodicity of $t_1\cdots t_{m-1}$
and the assumption $Z_i\neq m-1$ for $i\geq q$, given by the definition of the
index $q$, that $t_{Z_i}=t_i=t_{i\ {\rm mod}\ s}$, and therefore $Z_{i+1}=Z_i+1$
for all $i$, $q\leq i\leq l$. The periodicity of $t_1\cdots t_{m-1}$ also implies
$t_{l+1}=t_{Z_l+1}=t_Z+1$, which contradicts~\eqref{e3}.
\pfk

\subsection{Necessary condition for affine factor complexity of $u_\beta$}

We now show that if there exists a word $p$ which is both a proper
prefix and a proper suffix of $t_1\cdots t_{m-1}$, and $t_1\cdots
t_{m-1}$ is not an integer power of $p$, then the factor complexity of
$u_\beta$ is not an affine function. According to
Proposition~\ref{p:3} it suffices to find a left special factor
which is not a prefix of $u_\beta$.

For that, we use the fact that the words of ${\cal L}(u_\beta)$
code the patterns of $\Z_\beta^+$. Indeed, the word
$w=w_0w_1\cdots w_n\in\{0,1,\dots,m-1\}^*$ is a coding of the set 
$[x,y]\cap\Z_\beta$, if the distances between points of 
$[x,y]\cap\Z_\beta$ are consecutively $T_{\beta}^{w_0}(1)$, 
\dots, $T_{\beta}^{w_n}(1)$. With this, we can reformulate the 
main problem of this section in the language of $\Z_\beta^+$. 
Construction of a left special factor of $u_\beta$ which is not a 
prefix of $u_\beta$ is equivalent to the construction of 
$\beta$-integers $z$, $x_1$, $x_2$ such that

\begin{enumerate}
\item[(i)] the codings of the sets $[x_1,x_1+z]\cap \Z_\beta$, $[x_2,x_2+z]\cap
\Z_\beta$ and $[0,z]\cap \Z_\beta$ are equal to the same word $w$;
\item[(ii)] $x_1-{\rm pred}(x_1) \neq x_2-{\rm pred}(x_2)$;
\item[(iii)] $1={\rm succ}(x_1+z)-(x_1+z) = {\rm succ}(x_2+z)-(x_2+z)$;
\item[(iv)] $1\neq {\rm succ}(z)-z$.
\end{enumerate}

\noindent Note that as the distance $1=T_\beta^0(1)$ is coded by
the letter 0, conditions (i)-(iv) ensure that the word $w0\in{\cal
L}(u_\beta)$ is a left special factor of $u_\beta$ which is not a
prefix of $u_\beta$.

The construction of the suitable $\beta$-integers $z$,$x_1$,$x_2$
with the above properties, is the content of this section, we
shall however need some preparation.

Let $p=p_1\cdots p_s$, be a proper prefix and a proper suffix of
the word $t_1\cdots t_{m-1}$ of the minimal non-zero length. From
the Parry condition and the fact that $t_1\cdots t_{m-1}\neq p^k$
for $k\geq 2$ one can easily deduce that there exist words $p'$,
$q$, and a positive integer $r$ such that
\begin{equation}\label{cislo1}
d_\beta(1) = p^rp'qp1\,,
\end{equation}
where $p'$ is a prefix of $p$ and $|p|> |p'|:=j$, and $q$ is a
non-empty word starting with the letter $q_1<p_{j+1}$. Let us
mention that the words $p,p',q$ are words over the alphabet
$\{t_1,t_2,\dots,t_{m-1}\}$. Since $t_1\cdots t_{m-1}$ is not an
integer power of $p$, we must have
$$
pp'q\neq p'qp\,.
$$
As $|pp'q|=|p'qp|$, we can find a word $c\in\{t_1,t_2,\dots,
t_{m-1} \}^*$ and digits $h_1,h_2\in\{t_1,\dots,t_{m-1}\}$ such
that $h_1\neq h_2$, $h_1c$ is a suffix of $pp'q$, and $h_2c$ is a
suffix of $p'qp$. Note that since $q_1<p_{j+1}$, $c$ as a common
suffix of $pp'q$ and $p'qp$ must satisfy
\begin{equation}\label{eq:c}
|c|\leq |p|+|q|-1\,.
\end{equation}
Denote $h:=\min(h_1,h_2)$ and $A=|p^rp'q_1|=rs+j+1$. Then we can
define $\beta$-integers $x_1,x_2,z$ using their $\beta$-expansion
as
$$
\begin{array}{ccrcr}
z&:=&  p^rp'q_1\bullet  & + & hc0^A\bullet\\[2mm]
x_1&:=& p^{r}p'q0^A\bullet& -& hc0^A\bullet\\[2mm]
x_2&:=& p^rp'qp0^A\bullet& - & hc0^A\bullet
\end{array}
$$
Directly from the definition of $A$, $h$ and $c$, it follows that
the word $hcp^rp'q_1$ satisfies the Parry condition, and thus
$z=hcp^rp'q_1\bullet$ is a $\beta$-integer. In the same time, from
the definition of $h$ and $c$ it is obvious that subtraction in
the prescription for $x_1$ and $x_2$ can be performed digit-wise
and hence also $x_1,x_2\in\Z_\beta^+$.

We now prove that the above defined $z, x_1, x_2$ satisfy the
conditions (i) -- (iv).

\begin{trivlist}
\item[\underline{\bf ad(i)}]
In order to prove that the word $w$ coding the distances between
consecutive $\beta$-integers in the segment $[0,z]$, codes also
the segments $x_1,x_1+z$ and $x_2,x_2+z$, we use the following
lemma, which is a consequence of the fact that the infinite word
$u_\beta$ codes the distances between consecutive
$\beta$-integers.

\begin{lem}\label{3.20}
Let $x,z\in\Z_\beta^+$ such that
\begin{equation}\label{cislo2}
\hbox{for every } z'\in\Z_\beta^+,\ z'\leq z \hbox{ we have }
x+z'\in\Z_\beta^+\,.
\end{equation}
Then codings of $[0,z]\cap\Z_\beta$ and $[x,x+z]\cap\Z_\beta$
coincide.
\end{lem}

For $x=x_1$ we divide the verification of condition~\eqref{cislo2}
into three cases.

\begin{itemize}

\item
If $z'\in\Z_\beta$, $0\leq z'<hc0^A\bullet + p10^{A-s-1}\bullet$,
then the summation $x_1+z'$ can be performed digit-wise. The
result is again a string satisfying the Parry condition, and
therefore $x_1+z'\in\Z_\beta$.

\item
If $z'=hc0^A\bullet+p10^{A-s-1}\bullet$, then $ x_1+z' =
p^rp'q0^A\bullet + p10^{A-s-1}\bullet = 10^m0^{A-s-1}\bullet$.

\item
If $z'\in\Z_\beta$, $hc0^A\bullet + p10^{A-s-1}\bullet<z'\leq z$,
then $x_1+z'=10^m0^{A-s-1}\bullet + z''$, where $z''\in\Z_\beta$,
$0<z''\leq p^rp'q_1\bullet-p10^{A-s-1}\bullet$, and again by
digit-wise summation we obtain an admissible $\beta$-expansion of
$x_1+z'$.

\end{itemize}

In order to prove the condition~\eqref{cislo2} for $x=x_2$, we
again separate $z'\in\Z_\beta$, $z'\leq z$ into three cases. Now
the separating point is $z'=hc0^A\bullet + 10^{A-1}\bullet$. For
such $z'$ we have $x_2+z'=10^{m+A-1}\bullet$. The remaining cases
$x_2+z'$ can be solved by digit-wise summation, similarly as for
$x=x_1$.

\medskip
\item[\underline{\bf ad(ii)}]
For the proof of property (ii) we use another statement, which
allows one to determine the distance of an element of $\Z_\beta^+$
from its predecessor.

\begin{lem}\label{3.21}
Let the $\beta$-expansion of a $\beta$-integer $y$ be
$y_ny_{n-1}\cdots y_k0^k\bullet$, where $y_k\neq 0$, $k\in\N$.
Then $y-{\rm pred}(y)=T^{k'}_\beta(y)$, where
$k'\in\{0,1,\dots,m-1\}$ is such that $k'=k \!\mod m$.
\end{lem}

\pf
 If $k=0$, the statement is obvious. Assume that $k\geq 1$.
Denote $d_\beta^*(1)=\bigl(t_1t_2\cdots t_{m-1}0\bigr)^\omega$.
This is the lexicographically greatest word which is
lexicographically strictly smaller than $d_{\beta}(1)=t_1\cdots
t_{m-1}1$. Therefore the predecessor of $y$ has the
$\beta$-expansion of the form $y_n\cdots y_{k+1}(y_k-1)v\bullet$,
where $v$ is a prefix of $d_\beta^*(1)$ of length $|v|=k$. We thus
have
$$
y-{\rm pred}(y) \quad = \quad 10^k\bullet \ - \ v\bullet \quad =
\quad 0\bullet t_{k'+1}t_{k'+2}\cdots t_m  \quad = \quad
T^{k'}_\beta(1)\,,
$$
where $k'\in\{0,1,\dots,m-1\}$ is such that $k'\equiv k \!\mod m$.
The latter follows from the periodicity of $d_\beta^*(1)$ with
period $m$.
\pfk

The definition of $h$ implies that the number of 0's at the end of
$\beta$-expansions of $x_1,x_2$ differ modulo $m$. Therefore
property (ii) is valid.

\medskip
\item[\underline{\bf ad(iii)}] For verifying the property (iii) we have to show that
both $x_1+z+1$ and $x_2+z+1$ belong to $\Z_\beta$. Setting $d_\beta(1):=l$, we have 
$d_\beta(x_1+z)=10^l(p_1-1)p_2\cdots p_sp^{r-2}p'q_1$ if $r\geq2$, 
$d_\beta(x_1+z)=10^l(p'_1-1)p'_2\cdots p'_jq_1$ otherwise. In both cases, since 
$q_1<p_{j+1}$, the digit $1$ can be added at the last position of $d_\beta(x_1+z)$ 
without altering the validity of the Parry condition. The same argument holds for 
$d_\beta(x_2+z)=10^l(p_1-1)p_2\cdots p_sp^{r-1}p'q_1$.

\medskip
\item[\underline{\bf ad(iv)}] In order to prove that ${\rm succ}(z)\neq z+1$ we use
the statement which is a simple consequence of the proof of
Lemma~\ref{3.21}.

\begin{lem}
Let the $\beta$-expansion of a $\beta$-integer $y$ be
$y_ny_{n-1}\cdots y_0\bullet$. Denote by $k$ the maximal index
such that $y_{k-1}y_{k-2}\cdots y_{0}$ is a prefix of
$d_\beta^*(1)=\bigl(t_1t_2\cdots t_{m-1}0\bigr)^\omega$. Then
${\rm succ}(y)=y+T^{k'}_\beta(1)$, where $k'\in\{0,1,\dots,m-1\}$
is such that $k'\equiv k \!\mod m$.
\end{lem}

For $y=z=hcp^rp'q_1\bullet$, the index $k$ satisfies
$$
0<|p^rp'q_1|\leq k < |hcp^rp'q_1| \leq m\,,
$$
where we have used inequality~\eqref{eq:c}. Therefore
$T^k_\beta(1)\neq T^0_\beta(1)$.
\end{trivlist}

\section{Conclusions}

Among the words $u_\beta$ which have affine factor complexity are
words for which the R\'enyi expansion of unity in base $\beta$ is
of the form $d_\beta(1)=t_1t_2\cdots t_{m-1}1=p^k1$, for some
$k\geq 2$. If $p$ is a word of length 1, such words are
Arnoux-Rauzy, and thus have for each $n$ exactly one left special
and one right special factor of length $n$. If $p$ is of length
$|p|\geq 2$, then $u_\beta$ has for every $n\in\N$ one left
special and $|p|$ right special factors.

As a continuation of this paper, it would be interesting to study 
the factor complexity of a fixed point of a substitution defined by 
a non-simple Parry number.
It would also be interesting to compute explicitly the factor 
complexity in the non-affine case. In particular, is it possible 
that the factor complexity is ultimately affine, that is, 
${\cal C}(n)=an+b$ for $n\geq n_0$ ? Due to Lemma \ref{l:8}, there 
cannot exist finitely many maximal left special factors in the 
non-affine case, hence $a>m-1$ in such a case.

\section*{Acknowledgements}

The authors acknowledge financial support by Czech Science
Foundation GA \v{C}R 201/05/0169, by the grant LC06002 of the
Ministry of Education, Youth, and Sports of the Czech Republic.



\begin{thebibliography}{9}

\bibitem{cassaigne} J. Cassaigne.
Complexit\'e et facteurs sp\'eciaux. {\em Bull. Belg. Math. Soc.
Simon Stevin} {\bf 4} (1997), 67--88.


\bibitem{Durand} F. Durand.
Linearly recurrent subshifts have a finite number of non-periodic
subshift factors. {\em Ergodic Theory Dynam. Systems}, {\bf 20(4)}
(2000), 1061-1078.

\bibitem{fabre} S. Fabre.
Substitutions et $\beta$-syst\`emes de num\'eration. {\em Theoret.
Comput. Sci.} {\bf 137} (1995), 219--236.

\bibitem{Ferenczi} S. Ferenczi.
Complexity of sequences and dynamical systems. {\em Combinatorics
and number theory (Tiruchirappalli, 1996), Discrete Math.} {\bf
206 (1-3)} (1999), 145--154.

\bibitem{FrMaPe} Ch. Frougny, Z. Mas\'akov\'a, E. Pelantov\'a.
Complexity of infinite words associated with beta-expansions. {\em
RAIRO Theor. Inform. Appl.} {\bf 38} (2004), 163--185;
Corrigendum, {\em RAIRO Theor. Inform. Appl.} {\bf 38} (2004),
269--271.

\bibitem{lothaire} M. Lothaire.
{\it Algebraic combinatorics on words}. Cambridge University Press
(2002).

\bibitem{morse} Hedlund, G. A. and Morse, M.
Symbolic dynamics {II}. {S}turmian trajectories. {\it Amer. J.
Math.} {\bf 62} (1940), 1--42.

\bibitem{Pansiot} J-J. Pansiot.
Complexit\'e des facteurs des mots infinis engendr\'es par morphismes
it\'er\'es. {\it Automata, languages and programming (Antwerp, 1984)}.
Lecture Notes in Comput. Sci., Springer. {\bf 172} (1984),
380--389.

\bibitem{parry} W. Parry.
On the $\beta$-expansions of real numbers. {\em Acta Math. Acad.
Sci. Hungar.} {\bf 11} (1960), 401--416.

\bibitem{queffelec} M. Queff\'elec.
Substitution dynamical systems---spectral analysis. {\it Lecture
Notes in Mathematics, Springer-Verlag}. {\bf 1294} (1987),
xiv+240.

\bibitem{renyi} A. R\'enyi.
Representations for real numbers and their ergodic properties.
{\em Acta Math. Acad. Sci. Hungar.} {\bf 8} (1957), 477--493.

\bibitem{thurston} W.P. Thurston.
{\em Groups, tilings, and finite state automata}. Geometry
supercomputer project research report GCG1, University of
Minnesota (1989).


\end{thebibliography}
\end{document}